\documentclass[12pt]{amsart} 

\usepackage{amsmath} \usepackage{amssymb} \usepackage{amsthm}
\usepackage{amscd} \usepackage[all]{xy} 

\def\a{{\mathfrak{a}}} \def\b{{\mathfrak{b}}}
\def\F{{\mathbb{F}}} \def\J{{\mathcal{J}}} \def\m{{\mathfrak{m}}} \def\Z{{\mathbb{Z}}}\def\N{{\mathbb{N}}} \def\O{{\mathcal{O}}}
\def\Q{{\mathbb{Q}}} \def\R{{\mathbb{R}}} \def\C{{\mathbb{C}}}
\def\Hom{{\mathrm{Hom}}}  
  \def\Ann{{\mathrm{Ann}}}
\def\Spec{{\mathrm{Spec\; }}}

\theoremstyle{plain}
\newtheorem{thm}{Theorem}[section] 
\newtheorem{prop}[thm]{Proposition}

\newtheorem{propdef}[thm]{Proposition-Definition} 

\newtheorem{lem}[thm]{Lemma}
\theoremstyle{definition} 
\newtheorem{defn}[thm]{Definition}
\newtheorem{eg}[thm]{Example} 
\theoremstyle{remark}
\newtheorem{rem}[thm]{Remark}

\newtheorem{setup}[thm]{Set up}
\newtheorem*{cl}{Claim}

\newtheorem*{acknowledgement}{Acknowledgment}

\title{Formulas for multiplier ideals on singular varieties}
\author{Shunsuke Takagi}
\address{Faculty of Mathematics, Kyushu University, 
6-10-1 Hakozaki, Higashi-ku, Fukuoka, 812-8581 JAPAN}
\email{stakagi@math.kyushu-u.ac.jp}
\subjclass[2000]{13A35, 14B05}

\dedicatory{Dedicated to Professor~Kei-ichi~Watanabe on the~occasion of his~sixtieth~birthday.}

\begin{document}

\begin{abstract}
We prove a generalization of Demailly-Ein-Lazarsfeld's subadditivity formula \cite{DEL} and Musta\c t\v a's summation formula \cite{Mu3} for multiplier ideals to the case of singular varieties, using characteristic $p$ methods. 
As an application of our formula, we improve Hochster-Huneke's result \cite{HH4} on the growth of symbolic powers of ideals in singular affine algebras. 
\end{abstract}

\maketitle
\markboth{SHUNSUKE TAKAGI}{FORMULAS FOR MULTIPLIER IDEALS ON SINGULAR VARIETIES}

\section*{Introduction}
Multiplier ideals satisfy several nice properties related to vanishing theorems, making them an important tool in higher dimensional geometry. They are defined as follows: let $X$ be a $\Q$-Gorenstein normal variety over a field of characteristic zero and let $\a \subseteq \O_X$ be an ideal sheaf of $X$. Suppose that $f:\widetilde{X} \to X$ is a log resolution of $(X,V(\a))$, that is, $f$ is proper and birational, $\widetilde{X}$ is nonsingular and $f^{-1}V(\a)=F$ is a divisor with simple normal crossing support. If $K_{\widetilde{X}/X}$ is the relative canonical divisor of $f$, then the multiplier ideal of $\a$ with coefficient $t \in \R_{>0}$ is 
$$\J(\a^t)=\J(t \cdot \a)=f_*\O_{\widetilde{X}}(\lceil K_{\widetilde{X}/X}-tF \rceil) \subseteq \O_X.$$
The reader is referred to \cite{La} for general properties of multiplier ideals.

Demailly, Ein and Lazarsfeld \cite{DEL} proved a subadditivity property of multiplier ideals on smooth varieties, which states that the multiplier ideal of the product of two ideals is contained in the product of their individual multiplier ideals. 
Their strategy is as follows: first the problem is pulled back to the product $X \times X$ and then the desired formula on $X$ is recovered by restricting to the diagonal. This proof relies heavily on the fact that the diagonal embedding is a complete intersection and so it works only on nonsingular varieties. Indeed, several counterexamples to their formula are known on singular varieties (see \cite{TW}). 
Later Musta\c t\v a \cite{Mu3} also proved a formula for the multiplier ideal of a sum of two ideals analogous to Demailly-Ein-Lazarsfeld's formula.
Though his proof uses a more subtle argument, it follows essentially the same idea as above, and therefore the assumption that $X$ is nonsingular is also necessary for his proof. 
In this paper, we prove a generalization of the above two formulas for multiplier ideals to the case of singular varieties, using the theory of tight closure. 

The notion of tight closure is a powerful tool in commutative algebra introduced in $1980$'s by Hochster and Huneke \cite{HH1} using the Frobenius map. 
In this theory, the test ideal $\tau(R)$ plays a central role. 
Hara and Yoshida \cite{HY} introduced a generalization of the test ideal, the ideal $\tau(\a^t)$ associated to a given ideal $\a$ with real exponent $t> 0$, and they showed that this ideal corresponds to the multiplier ideal $\J(\a^t)$ via reduction to characteristic $p \gg 0$. 
By virtue of their result, we can study local properties of the multiplier ideal $\J(\a^{t})$ by examining those of the ideal $\tau(\a^{t})$. 
First we prove subadditivity and summation formulas involving the Jacobian ideal for $\tau(\a^t)$  on singular affine algebras, by taking advantage of arguments developed in \cite{HH4} and \cite{HT}. 
After that, via the above correspondence, we obtain the following formulas for multiplier ideals. 

\begin{thm}
Let $X$ be a normal $\Q$-Gorenstein variety over a field $K$ of characteristic zero and let $\mathfrak{J}(X/K)$ be the Jacobian ideal sheaf of $X$ over $K$. 
Let $\a, \b \subseteq \O_X$ be two nonzero ideal sheaves on $X$ and fix any real numbers $s, t>0$.  
\begin{enumerate}
\item \textup{(Subadditivity formula)}
$$\mathfrak{J}(X/K)\J(\a^t\b^s) \subseteq \J(\a^t)\J(\b^s).$$
\item \textup{(Summation formula)}
$$\J((\a+\b)^t)=\sum_{\lambda+\mu=t}\J(\a^{\lambda}\b^{\mu}).$$ 
In particular, 
$$\mathfrak{J}(X/K) \J((\a+\b)^t) \subseteq \sum_{\lambda+\mu=t}\J(\a^{\lambda}) \J(\b^{\mu}).$$
\end{enumerate}
\end{thm}

As an application of our subadditivity formula, we prove a generalization of the result due to Ein-Lazarsfeld-Smith \cite{ELS} and Hochster-Huneke \cite{HH4} concerning the growth of symbolic powers of ideals in affine $K$-algebras. Hochster and Huneke proved that if $h$ is the largest height of an associated prime of an ideal $\a \subseteq R$ and $J$ is the Jacobian ideal of $R$ over $K$, then $J^{n+1}\a^{(hn)} \subseteq \a^n$ for all $n \ge 0$. 
Here, if $W$ is the complement of the union of the associated primes of $\a$, $\a^{(t)}$ denotes the contraction of $\a^tR_W$ to $R$, where $R_W$ is the localization of $R$ at the multiplicative system $W$.
They then asked whether the exponent $n+1$ used on the Jacobian ideal is best possible. 
Recently Hara \cite{Ha} defined a variant of $\tau(\a^t)$ corresponding to asymptotic multiplier ideals and applied it to give an alternative proof of their result in the case where $R$ is regular.  
Employing the same strategy (which can be traced back to the method in \cite{ELS}), we use our subadditivity formula to improve Hochster-Huneke's result by decreasing the exponent on the Jacobian ideal by $1$. 
\begin{thm}
Let $R$ be an equidimensional reduced affine algebra over a perfect field $K$ of positive characteristic and $J=\mathfrak{J}(R/K)$ be the Jacobian ideal of $R$ over $K$. 
Let $\a \subseteq R$ be any ideal which is not contained in any minimal prime ideal. 
Let $h$ be the largest analytic spread of $\a R_P$ as $P$ runs through the associated primes of $\a$. 
Then, for every integer $m \ge 0$ and every integer $n \ge 1$, 
$$J^n\a^{(hn+mn)} \subseteq (\a^{(m+1)})^n. $$
\end{thm}

\section{Preliminaries on a generalization of test ideals}
In this section, we briefly review the definition and fundamental properties of a generalization of test ideal $\tau(\a^t)$ which we need later. We refer the reader to \cite{HY} and \cite{HT} for the proofs. 

In this paper, all rings are excellent reduced commutative rings with unity. 
For a ring $R$, we denote by $R^{\circ}$ the set of elements of $R$ which are not in any minimal prime ideal. 
Let $R$ be a ring of characteristic $p >0$ and $F\colon R \to R$ the Frobenius map which sends $x \in R$ to $x^p \in R$. 
For an integer $e > 0$, the ring $R$ viewed as an $R$-module via the $e$-times iterated Frobenius map $F^e \colon R \to R$ is denoted by ${}^e\! R$. 
Since $R$ is assumed to be reduced, we can identify $F^e \colon R \to {}^e\! R$ with the natural inclusion map $R \hookrightarrow R^{1/p^e}$. We say that $R$ is {\it F-finite} if ${}^1
\! R$ (or $R^{1/p}$) is a finitely generated $R$-module. 

Let $R$ be a ring of characteristic $p > 0$ and $M$ an $R$-module. 
For each integer $e > 0$, we denote $\F^e(M) = \F_R^e(M) := M \otimes_R {}^e\! R$ and regard it as an $R$-module by the action of $R$ on ${}^e\! R$ from the right. 
Then we have the induced $e$-times iterated Frobenius map $F_M^e \colon M \to \F^e(M)$. 
The image of $z \in M$ via this map is denoted by $z^q:= F^e(z) 
\in \F^e(M)$, where $q=p^e$. For an $R$-submodule $N$ of $M$, we denote by $N^{[q]}_M$ the 
image of the induced map $\F^e(N) \to \F^e(M)$. 

Now we recall the definition of $\a^{t}$-tight closure.
\begin{defn}[\textup{\cite[Definition 6.1]{HY}}]\label{ta}
Let $\a$ be an ideal of a ring $R$ of characteristic $p>0$ such that $\a \cap R^{\circ} \ne \emptyset$ and let $N \subseteq M$ be $R$-modules. 
Given a real number $t> 0$, the {\it $\a^{t}$-tight closure} $N^{*\a^{t}}_M$ of $N$ in $M$ is defined to be the submodule of $M$ consisting of all elements $z \in M$ for which there exists $c \in R^{\circ}$ such that 
$$cz^q\a^{\lceil tq \rceil} \subseteq N^{[q]}_M$$ 
for all large $q = p^e$, where $\lceil tq \rceil$ is the least integer which is greater than or equal to $tq$. 
The $\a^{t}$-tight closure of an ideal $I \subseteq R$ is just defined by $I^{*\a^{t}} = I^{*\a^{t} }_R$.
\end{defn}

\begin{rem}
(1) When $\a$ is the unit ideal, $\a^t$-tight closure is nothing but classical tight closure. 
The reader is referred to Hochster-Huneke's original paper \cite{HH1} for the classical tight closure theory. 

(2) (Several exponents) 
The real exponent $t$ for $\a^t$-tight closure in Definition 1.1 is just a formal notation, but it is compatible with ``real" powers of the ideal. 
Namely, if $\b = \a^n$ for $n \in \N$, then $\a^t$-tight closure is the same as $\b^{t/n}$-tight closure. This allows us to extend the definition to several real exponents:
given ideals $\a_1, \dots, \a_r \subseteq R$ with $\a_i \cap R \ne \emptyset$ and real numbers $t_1, \dots, t_r > 0$, 
if $t_i=tn_i$ for $t \in \R_{>0}$ and $n_i \in \N$ with $i=1, \dots, r$, we can define $\a_1^{t_1} \dots \a_r^{t_r}$-tight closure to be $(\a_1^{n_1} \dots \a_r^{n_r})^t$-tight closure. If $N$ is a submodule of $R$-module $M$, then an element $z \in M$ is in the $\a_1^{t_1} \dots \a_r^{t_r}$-tight closure $N^{*\a_1^{t_1} \dots \a_r^{t_r}}$ of $N$ in $M$ if and only if there exists $c \in R^{\circ}$ such that $c z^q \a_1^{\lceil t_1 q \rceil} \dots \a_r^{\lceil t_r q \rceil} \subseteq N_M^{[q]}$ for all $q=p^e \gg 0$. 
\end{rem}

Now using the $\a^{t}$-tight closure of the zero submodule, we define two ideals $\tau(\a^t)$ and $\widetilde{\tau}(\a^t)$.
\begin{propdef}[\textup{\cite[Proposition-Definition 1.9]{HY}, cf.\cite{HH1}}]\label{taudef}
Let $R$ be an excellent reduced ring of characteristic $p > 0$, let $\a \subseteq R$ be an ideal with $\a \cap R^{\circ} \ne \emptyset$  and let $t > 0$ be a real number. 
Let $E =\bigoplus_{\m} E_R(R/\m)$ be the direct sum, taken over all maximal ideals $\m$ of $R$, of the injective hulls of the residue fields $R/\m$. 
Then the following ideals are equal to each other 
and we denote them by $\tau(\a^{t})$. 
\begin{enumerate}
\renewcommand{\labelenumi}{\textup{(\roman{enumi})}}
\item $\displaystyle\bigcap_M \Ann_R(0^{*\a^{t}}_M)$, where $M$ runs through all finitely generated $R$-modules. 

\item $\displaystyle\bigcap_{M\subset E} \Ann_R(0^{*\a^{t}}_M)$, where $M$ runs through all finitely generated submodules of $E$. 
\item $\displaystyle\bigcap_{J\subseteq R} (J:J^{*\a^{t}})$, where $J$ runs through all ideals of $R$. 
\end{enumerate}
\end{propdef}

\begin{defn}[\textup{\cite[Definition 1.4]{HT}}]
Let $E=\bigoplus_{\m}E_R(R/\m)$ be as in Definition-Theorem \ref{taudef}. 
Then we define the ideal $\widetilde{\tau}(\a^t)$ by 
$$\widetilde{\tau}(\a^t)=\Ann_R(0_E^{*\a^t}).$$
\end{defn}

In general $\widetilde{\tau}(\a^t) \subseteq \tau(\a^t)$, and we have no example in which these two ideals disagree. 
\begin{thm}[\textup{\cite[Theorem 1.13]{HY}, cf.\cite{AM}}]
Let $R$ be an excellent normal $\Q$-Gorenstein ring of characteristic $p>0$. 
Then, for every real number $t>0$ and every ideal $\a \subseteq R$ such that $\a \cap R^{\circ} \ne \emptyset$,
$$\tau(\a^t)=\widetilde{\tau}(\a^t).$$
\end{thm}

Next we define the notion of $\a^t$-test elements, which is useful to study the ideal $\tau(\a^t)$. We refer the reader to \cite[Theorem 1.7, 6.4]{HY} and \cite[Corollary 2.4]{HT} for the existence of $\a^t$-test elements. 

\begin{defn}[\textup{\cite[Definition 6.3]{HY}}]\label{testdef}
Let $\a$ be an ideal of a ring $R$ of characteristic $p > 0$ such that $\a \cap R^{\circ} \ne \emptyset$ and let $t> 0$ be a real number. 
An element $d \in R^{\circ}$ is called an {\it $\a^{t}$-test element} if for every finitely generated $R$-module $M$ and $z \in M$, the following holds: $z \in 0_M^{*\a^{t}}$ if and only if $dz^q\a^{\lceil tq \rceil}= 0$ for all powers $q = p^e$ of $p$.
\end{defn}

A criterion for the existence of $\a^t$-test elements given in \cite[Corollary 2.4]{HT} also works for ``test elements for $0_E^{*\a^t}$.'' 

\begin{thm}\label{exist}
Let $(R,\m)$ be an F-finite reduced local ring of characteristic $p>0$ and let $E=E_R(R/\m)$ be the injective hull of the residue field $R/\m$. 
\begin{enumerate}
\item \textup{(cf. \cite[Theorem 3.3]{HH2})}
Let $d \in R^{\circ}$ be an element such that the localized ring $R_d$ is regular. 
Then some power $d^n$ of $d$ belongs to $\widetilde{\tau}(R)$. 

\item \textup{(cf. \cite[Corollary 2.4]{HT})}
Let $d \in R^{\circ}$ be an arbitrary element of $\widetilde{\tau}(R)$. 
Then, for all ideals $\a \subseteq R$ such that $\a \cap R^{\circ} \ne \emptyset$ and for all real numbers $t> 0$, the following holds: for every element $z \in E$, $z \in 0_E^{*\a^{t}}$ if and only if $dz^q\a^{\lceil tq \rceil}= 0$ for all powers $q = p^e$ of $p$.
\end{enumerate}
\end{thm}
\begin{proof}
The proofs are essentially the same as those of \cite[Theorem 3.3]{HH2} and \cite[Corollary 2.4]{HT}. 
\end{proof}

The following lemma gives a characterization of the ideal $\widetilde{\tau}(\a^t)$. 
\begin{lem}[\textup{\cite[Lemma 2.1]{HT}}]\label{key lemma}
Let $(R,\m)$ be an F-finite reduced local ring of characteristic $p > 0$, 
let $\a \subseteq R$ be an ideal such that $\a \cap R^{\circ} \ne \emptyset$ and let $t > 0$ be a real number. 
Let $d \in R^{\circ}$ be an element of $\widetilde{\tau}(R)$. 
Fix a system of generators $x^{(e)}_1,\ldots, x^{(e)}_{r_{e}}$ of $\a^{\lceil tq \rceil}$ for each $q = p^e$. Then an element $c \in R$ belongs to $\widetilde{\tau}(\a^t)$ if and only if there exist an integer $e'>0$ and $R$-linear maps $\phi^{(e)}_i \in \Hom_R(R^{1/p^e},R)$ for $0 \le e \le e'$ and $1 \le i \le r_e$ such that 
$$c=\sum_{e=0}^{e'}\sum_{i=1}^{r_e} \phi^{(e)}_i((dx_i^{(e)})^{1/p^e}).$$
\end{lem}

As applications of Lemma \ref{key lemma}, we can show that the ideal $\widetilde{\tau}(\a^t)$ commutes with localization and completion. 
\begin{prop}[\textup{\cite[Proposition 3.1]{HT}}]\label{loc}
Let $(R,\m)$ be an F-finite reduced local ring of characteristic $p >0$, let $\a \subseteq R$ be an ideal such that $\a \cap R^{\circ} \ne \emptyset$ and let $t > 0$ be a real number. 
Let $W$ be a multiplicatively closed subset of $R$. Then
$$\widetilde{\tau}((\a R_W)^t) = \widetilde{\tau}(\a^t)R_W.$$
\end{prop}

\begin{prop}[\textup{\cite[Proposition 3.2]{HT}}]\label{com}
Let $(R,\m)$ be an F-finite reduced local ring of characteristic $p>0$, 
let $\a \subseteq R$ be an ideal such that $\a \cap R^{\circ} \ne \emptyset$ and let $t > 0$ be a real number. 
Let $\widehat{R}$ denote the $\m$-adic completion of $R$. Then
$$\widetilde{\tau}((\a\widehat{R})^t)=\widetilde{\tau}(\a^t)\widehat{R}.$$
\end{prop}

Also we have Skoda's theorem for $\widetilde{\tau}(\a^t)$ with the aid of Lemma \ref{key lemma}. 
\begin{thm}[\textup{\cite[Theorem 4.2]{HT}}]\label{skoda}
Let $R$ be an F-finite reduced ring of characteristic $p > 0$.  
Let $\a \subseteq R$ be an ideal such that $\a \cap R^{\circ} \ne \emptyset$ and assume that $\a$ has a reduction generated by $l$ elements. 
Let $\b$ be an ideal of $R$ such that $\b \cap R^{\circ} \ne \emptyset$ and let $t>0$ be a real number. 
Then 
$$\widetilde{\tau}(\a^l\b^t) = \widetilde{\tau}(\a^{l-1}\b^t)\a.$$
\end{thm}

\begin{setup}\label{reduction}
Let $R$ be an algebra essentially of finite type over a field $k$ of characteristic zero, let $\a \subseteq  R$ be an ideal and let $t>0$ be a real number. One can choose a finitely generated $\Z$-subalgebra $A$ of $k$ and a subalgebra $R_A$ of $R$ essentially of finite type over $A$ such that the natural map $R_A \otimes_A k \to R$ is an isomorphism and ${\a_A}R=\a$ where ${\a_A}=\a \cap R_A \subseteq R_A$. 
Given a closed point $s \in \Spec A$ with residue field $\kappa=\kappa(s)$, we denote the corresponding fibers over $s$ by $R_{\kappa}, {\a_{\kappa}}$. 
Then we refer to such $(\kappa, R_{\kappa}, {\a_{\kappa}})$ for a general closed point $s \in \Spec A$ with residue field $\kappa=\kappa(s)$ of sufficiently large characteristic $p \gg 0$ as ``{\it reduction to characteristic $p \gg 0$}'' of $(k,R, \a)$, and the pair $(R_{\kappa}, {\a_{\kappa}}^{t})$ inherits the properties possessed by the original pair $(R, \a^{t})$.
Furthermore, given a log resolution $f:\widetilde{X} \to X=\Spec R$ of $(X, \a^{t})$, we can reduce this entire setup to characteristic $p \gg 0$. 
\end{setup}
Hara and Yoshida prove the multiplier ideal $\J(\a^t)$ coincides, after reduction to characteristic $p \gg 0$, with the ideal $\tau(\a^t)$ (or equivalently $\widetilde{\tau}(\a^t)$). 
\begin{thm}[\textup{\cite[Theorem 6.8]{HY}, cf.\cite{Ha2}, \cite{Sm}}]\label{multiplier}
Let $(R,\m)$ be a $\Q$-Gorenstein normal local ring essentially of finite type over a field of characteristic zero. 
Let $\a \subseteq R$ be a nonzero ideal and let $t>0$ be a real number. 
Then, after reduction to characteristic $ p \gg 0$, 
$$\tau(\a^{t})=\widetilde{\tau}(\a^t)=\J(\a^{t}).$$
\end{thm}

\section{Subadditivity formula on singular varieties}
In this section, we consider a generalization of Demailly-Ein-Lazarsfeld's subadditivity formula \cite{DEL} for multiplier ideals to the case of singular varieties via the ideal $\tau(\a^t)$. 
The following two propositions are key propositions. 
\begin{prop}\label{complete}
Let $(R,\m)$ be a complete reduced local ring of characteristic $p>0$, let $\a \subseteq R$ be an ideal of positive height and let $t>0$ be a real number. 
Then,  for all ideals $\b \subseteq R$ of positive height and for all real numbers $s>0$, 
$$\tau(\a^t)^{*\a^t}\tau(\a^t\b^s) \subseteq \tau(\a^t)\tau(\b^s).$$
\end{prop}
\begin{proof}
The proof is similar to that of \cite[Theorem 4.1]{HT}. 
First we will see that 
$$(0_M^{\a^t \b^s}: \tau(\a^t)^{*\a^t})_M \supseteq (0_M^{\b^s}:\tau(\a^t))_M $$
for any (not necessarily finitely generated) $R$-module $M$. 
Let $z \in (0_M^{\b^s}:\tau(\a^t))_M$, that is, $\tau(\a^t)z \in 0_M^{\b^s}$. 
Then there exists $c \in R^{\circ}$ such that $c\b^{\lceil sq \rceil} \tau(\a^t)^{[q]}z^q=0$ in $\F^e(M)$ for all $q=p^e \gg 0$. 
On the other hand, by definition, there exists $d \in R^{\circ}$ such that $d\a^{\lceil tq \rceil}(\tau(\a^t)^{*\a^t})^{[q]} \subseteq \tau(\a^t)^{[q]}$ for all $q=p^e \gg 0$. 
Hence, one has  
$$cd\a^{\lceil tq \rceil} \b^{\lceil sq \rceil}(\tau(\a^t)^{*\a^t})^{[q]}z^q \subseteq  c\b^{\lceil sq \rceil}\tau(\a^t)^{[q]}z^q=0$$
in $\F^e(M)$ for all $q=p^e \gg 0$, namely $z \in (0_M^{\a^t \b^s}: \tau(\a^t)^{*\a^t})_M$. 
Thus we have $(0_M^{\a^t \b^s}: \tau(\a^t)^{*\a^t})_M \supseteq (0_M^{\b^s}:\tau(\a^t))_M $. 

Now assume that $(R,\m)$ is a complete local ring and let $E=E_R(R/\m)$ be the injective hull of the residue field $R/\m$. 
Then by Matlis duality, $\Ann_E(\tau(\b^s))$ is equal to the union of $0^{*\b^s}_M$ taken over all finitely generated $R$-submodules $M$ of $E$. 
Hence, if $z \in \Ann_E(\tau(\a^t)\tau(\b^s))$, then there exists a finitely generated submodule $M \subset E$ such that $z \in (0^{*\b^s}_M \colon \tau(\a^t))_E$. 
Replacing $M$ by $M+Rz \subset E$, one has $z \in (0^{*\b^s}_M \colon \tau(\a^t))_M \subseteq (0^{*\a^t\b^s}_M:\tau(\a^t)^{*\a^t})_M$. 
Consequently, $\Ann_E(\tau(\a^t)\tau(\b^s))$ is contained in the union of $(0^{*\a^t\b^s}_M:\tau(\a^t)^{*\a^t})_M$ taken over all finitely generated submodules $M \subset E$. Therefore, by Matlis duality again, 
\begin{align*}
\tau(\a^t)^{*\a^t}\tau(\a^t\b^s) 
& = \bigcap_{M\subset E} \Ann_R ((0^{*\a^t\b^s}_M:\tau(\a^t)^{*\a^t})_M)\\
& \subseteq \Ann_R(\Ann_E(\tau(\a^t)\tau(\b^s)))\\
& = \tau(\a^t)\tau(\b^s).
\end{align*}
\end{proof}

\begin{prop}\label{F-finite}
Let $R$ be an F-finite reduced ring of positive prime characteristic $p$. 
Let $\a$ be an ideal of $R$ such that $\a \cap R^{\circ} \ne \emptyset$ and let $t>0$ be a real number.  Then,  for all ideals $\b \subseteq R$ with $\b \cap R^{\circ} \ne \emptyset$ and  for all real numbers $s>0$, 
$$\widetilde{\tau}(\a^t)^{*\a^t} \widetilde{\tau}(\a^t\b^s) \subseteq \widetilde{\tau}(\a^t)\widetilde{\tau}(\b^s). $$
\end{prop}

\begin{proof}
Since taking annihilator is preserved under localization, we may assume that $R$ is a local ring with the maximal ideal $\m$. 
By Proposition \ref{com}, we can also assume that $R$ is complete. 
By the proof of Proposition \ref{complete}, one has $(0_E^{\a^t \b^s}: \widetilde{\tau}(\a^t)^{*\a^t})_E \supseteq (0_E^{\b^s}:\widetilde{\tau}(\a^t))_E$, where $E=E_R(R/\m)$ is the injective hull of the residue field $R/\m$. 
By Matlis duality, we have 
\begin{align*}
\widetilde{\tau}(\a^t)^{*\a^t} \widetilde{\tau}(\a^t\b^s)=\Ann_R((0_E^{\a^t \b^s}: \widetilde{\tau}(\a^t)^{*\a^t})_E) &\subseteq \Ann_R((0_E^{\b^s}:\widetilde{\tau}(\a^t))_E)\\
&=\widetilde{\tau}(\a^t)\widetilde{\tau}(\b^s). 
\end{align*}
\end{proof}

\begin{rem}
The proof of Proposition \ref{complete} (resp. Proposition \ref{F-finite}) tells us that  $I^{*\a^t} {\tau}(\a^t\b^s) \subseteq I \tau(\b^s)$ (resp. $I^{*\a^t} \widetilde{\tau}(\a^t\b^s) \subseteq I \widetilde{\tau}(\b^s)$) holds for any ideal $I \subseteq R$ under the assumption of Proposition \ref{complete} (resp. Proposition \ref{F-finite}). 
\end{rem}

As one of applications of Propositions \ref{complete} and \ref{F-finite},  we give a new proof of Hara-Yoshida's subadditivity formula \cite[Thoerem 4.5]{HY} for $\tau(\a^t)$ on regular rings. 
\begin{thm}[\textup{\cite[Theorem 4.5]{HY}}]
Let $R$ be a complete regular local ring of characteristic $p>0$ or F-finite regular ring of characteristic $p>0$. 
Then for any two ideals $\a$, $\b$ of $R$ and for any two positive real numbers $t,s$, 
$$\tau(\a^t\b^s) \subseteq \tau(\a^t)\tau(\b^s).$$ 
\end{thm}

\begin{proof}
Thanks to Propositions \ref{complete} and \ref{F-finite}, it suffices to show that $\tau(\a^t)^{*\a^t}$ is the unit ideal for every ideal $\a \subseteq R$ and every real number $t>0$. 
We may assume without loss of generality that $(R,\m)$ is a complete regular local ring as in the proof of Proposition \ref{F-finite}. 
Let $E$ be the injective hull of the residue field of $R$. 
Since $R$ is regular, $E \cong H^d_{\m}(R)$ and $E \otimes_R {}^e R \cong H^d_{\m}({}^e R)$.
We can identify $E$ with $E \otimes_R {}^e R$ via the identification of $R$ with ${}^e R$. 
Since the $e$-times iterated Frobenius map $F^e:R \to {}^e R$ is flat (because $R$ is regular: cf. \cite{Ku}), via the identification of $R \cong {}^e R$, we have $(0:\tau(\a^t))_E \otimes_R {}^e R \cong (0:\tau(\a^t)^{[q]})_E$ in $E \otimes_R {}^e R \cong E$. 
Let $F^e_E:E \to E \otimes_R {}^e R \cong E$ be the $e$-times iterated Frobenius map induced on $E$. Then $F^e_E((0:\tau(\a^t))_E)$ generates $(0:\tau(\a^t)^{[q]})_E$ in $E$. 

On the other hand, $(0:\tau(\a^t))=0_E^{\a^t}$ by Matlis duality, because $R$ is complete. 
Since the unit $1$ is an $\a^t$-test element by \cite[Theorem 1.7]{HY}, 
 $\a^{\lceil tq \rceil}F^e_E((0:\tau(\a^t))_E)=0$ in $E$ for all $q=p^e>0$. 
Thus $\a^{\lceil tq \rceil}(0:\tau(\a^t)^{[q]})_E=0$ for all $q=p^e>0$, which implies $\a^{\lceil tq \rceil} \subseteq \tau(\a^t)^{[q]}$ for all $q=p^e>0$ by Matlis duality again. 
\end{proof}

Next we consider a subadditivity property on singular affine algebras. 
The following result due to Hochster and Huneke is very useful to handle the Jacobian ideal. 
\begin{thm}[\textup{\cite[Theorem 3.4]{HH4}, cf.\cite{HH5}}]\label{jacobianlem}
Let $R$ be a geometrically reduced equidimensional affine algebra over a field $K$ of characteristic $p>0$. 
Let $L$ be an infinite extension field of $K$, and let $R_L=R \otimes_K L$. 
Let $\mathfrak{J}(R_L/L)$ be the Jacobian ideal of $R_L$ over $L$.
Then, for any element $c \in \mathfrak{J}(R_L/L)$, there exists a regular subring $A$ of $R_L$ (depending on $c$) such that $R_L$ is module-finite and generically \'etale over $A$.
Moreover, $A$ satisfies the following properties:
\begin{enumerate}
\item $c(R_L)^{1/q} \subseteq A^{1/q}[R_L]$ for every $q=p^e>0$. \\
\item $A^{1/q}[R_L] \cong A^{1/q} \otimes_A R_L$ is flat over $R_L$ for every $q=p^e>0$.
\end{enumerate}
\end{thm}

\begin{lem}\label{jacobian}
Let $R$ be an equidimensional reduced affine algebra over a perfect field $K$ of characteristic $p>0$. 
Let $\mathfrak{J}(R/K)$ be the Jacobian ideal of $R$ over $K$. 
Then, for all ideals $\a \subseteq R$ with $\a \cap R^{\circ} \ne \emptyset$ and  for all real numbers $t>0$, the ideal $\mathfrak{J}(R/K)$ is contained in the $\a^t$-tight closure $\widetilde{\tau}(\a^t)^{*\a^t}$ of $\widetilde{\tau}(\a^t)$. 
\end{lem}

\begin{proof}
Let $x$ be an indeterminate over $K$, let $L=K(x)$ and let $R_L=R \otimes_K L$. 
Since $K$ is a perfect field, $R$ and $R_L$ are F-finite. 
Fix an arbitrary element $c \in \mathfrak{J}(R/K)$ and any power $q$ of $p$. Then, by Theorem \ref{jacobianlem}, there exists a flat $R_L$-algebra $S_q$ such that $R_L \subseteq S_q \subseteq (R_L)^{1/q}$ and $c(R_L)^{1/q} \subseteq S_q$. 
First note that the issues are unaffected by replacing $R$ by $R_L$. 
In order to check this, we need the following claim. 

\begin{cl}[\textup{cf.\cite{HH3}}]
$$\widetilde{\tau}((\a R_L)^t)^{[q]}\cap R = \widetilde{\tau}(\a^t)^{[q]}. $$
\end{cl}

\begin{proof}[Proof of Claim]
Fix any maximal ideal $\m$ of $R$. 
Since $K$ is algebraically closed in $L$, $R_L/\m R_L=R/\m \otimes_K L$ is a field, that is, $\m R_L$ is a maximal ideal of $R_L$. 
Let $S_L=(R_L)_{\m R_L}$ and let $S=R_\m$. 
We can check the assertion locally, so it is enough to show that $\widetilde{\tau}((\a S_L)^t) = \widetilde{\tau}((\a S)^t)S_L$, because $S_L$ is faithfully flat over $S$. 
Let $E_L$ (resp. $E$) be the injective hull of the residue field of $S_L$ (resp. $S$). 
Then $E_L = E \otimes _K L $. Hence it is easy to see that  $0^{*(\a S_L)^t}_{E_L} = 0^{*(\a S)^t}_E \otimes_K L$ and we obtain the assertion. 
\end{proof}

By the above claim, the condition $\mathfrak{J}(R/K) \subseteq \widetilde{\tau}(\a^t)^{*\a^t}$ is equivalent to saying that 
$\mathfrak{J}(R/K)R_L = \mathfrak{J}(R_L/L) \subseteq \widetilde{\tau}((\a R_L)^t)^{*(\a R_L)^t}$,  because $R_L$ is faithfully flat over $R$.  
Thus we may assume without loss of generality that there exists a flat $R$-algebra $S_q$ such that $R \subseteq S_q \subseteq R^{1/q}$ and $c R^{1/q} \subseteq S_q$. 
Since $R$ is F-finite, $S_q$ is finite free as $R$-module.
Moreover, we can check the assertion locally, even passing to completion by virtue of Proposition \ref{com}. Thus we may reduce to the case where $R$ is a complete local ring with the maximal ideal $\m$. 

Let $E=E_R(R/\m)$ be the injective hull of the residue field $R/\m$. 
By Matlis duality, $(0:\widetilde{\tau}(\a^t))_{E}=0^{\a^t}_{E}$,
because $R$ is a complete local ring. 
Let $F^e_{E}:E \to E \otimes_{E} R^{1/q}$ be the $e$-times iterated  Frobenius map induced on $E$. 
Let $d \in R^{\circ}$ be an element of $\widetilde{\tau}(R)$ (such an element exists by Theorem \ref{exist} (1)). 
Then, by definition and Theorem \ref{exist} (2),  one has $d^{1/q}\a^{\lceil tq \rceil/q}F^e_{E}((0:\widetilde{\tau}(\a^t))_{E})=0$ in $E
\otimes_{E} R^{1/q}$. 
Here we consider the following $R$-module homomorphism.
$$\phi_q: E \xrightarrow{F^e_{E}}  E \otimes_{R} R^{1/q} \xrightarrow{\times c} E \otimes_{R}  S_q.$$
Since $\phi_q$ factors through $F^e_{E}$,   one has $d^{1/q} \a^{\lceil tq \rceil /q} \phi_q ((0:\widetilde{\tau}(\a^t))_{E})=c(d^{1/q}\a^{\lceil tq \rceil/q}F^e_{E}((0:\widetilde{\tau}(\a^t))_{E}))=0$ in $E \otimes_{R}  S_q$. 
On the other hand, since $S_q$ is flat over $R$,  
$$(0:\widetilde{\tau} (\a^t))_{E} \otimes_{R} S_q = (0 : \widetilde{\tau} (\a^t) S_q)_{E \otimes_{R} S_q}. $$
Hence $\phi_q ((0:\widetilde{\tau}(\a^t))_{E})$ generates $c (0 : \widetilde{\tau} (\a^t)S_q)_{E \otimes_{R} S_q}$ in $E \otimes_{R}  S_q$, because $\phi_q$ is factorized into $E \to  E \otimes_{R} S_q \xrightarrow{\times c} E \otimes_{R}  S_q$. 
Thus 
$$d^{1/q} \a^{\lceil tq \rceil /q} c (0 : \widetilde{\tau} (\a^t)S_q)_{E \otimes_{R} S_q}=0$$ in 
$E \otimes_{R}  S_q$. 
Now we have an isomorphism $\Hom_{R}(S_q, E) \cong E \otimes_{R} S_q$, because $S_q$ is a finite free $R$-module.  
We apply Matlis duality to $S_q$ via this isomorphism so that 
$$d^{1/q} \a^{\lceil tq \rceil /q} c \subseteq \Ann_{S_q}(0 : \widetilde{\tau} (\a^t)S_q)_{E \otimes_{R} S_q}=\widetilde{\tau} (\a^t)S_q \subseteq \widetilde{\tau} (\a^t)R^{1/q},$$ 
namely $d\a^{\lceil tq \rceil}c^q \in \widetilde{\tau}(\a^t)^{[q]}$. 
\end{proof}

As a direct consequence of Lemma \ref{jacobian}, we obtain a subadditivity formula involving the Jacobian ideal for $\widetilde{\tau}(\a^t)$ on singular affine algebras. 
\begin{thm}\label{subtau}
Let $R$ be an equidimensional reduced affine algebra over a perfect field $K$ of characteristic $p>0$. 
Let $\mathfrak{J}(R/K)$ be the Jacobian ideal of $R$ over $K$. 
Then, for any two ideals $\a, \b$ of $R$ such that $\a \cap R^{\circ} \ne \emptyset$ and  $\b \cap R^{\circ} \ne \emptyset$, and for any two positive real numbers $s, t$, 
$$\mathfrak{J}(R/K) \widetilde{\tau}(\a^t\b^s) \subseteq  \widetilde{\tau}(\a^t)\widetilde{\tau}(\b^s). $$
\end{thm}

\begin{proof}
Apply Lemma \ref{jacobian} to Proposition \ref{F-finite}. 
\end{proof}

\begin{rem}
Let $R$ be as in Theorem \ref{subtau} and let $I=\bigcap_{\a \subseteq R}\bigcap_{t>0} \widetilde{\tau}(\a^t)^{*\a^t} \subseteq R$, where $\a$ runs through all ideals of $R$ and $t$ runs through all positive real numbers. Then, by Proposition \ref{F-finite}, $I \widetilde{\tau}(\a^t\b^s) \subseteq  \widetilde{\tau}(\a^t)\widetilde{\tau}(\b^s)$ holds for any two ideals $\a, \b$ of $R$ such that $\a \cap R^{\circ} \ne \emptyset$ and  $\b \cap R^{\circ} \ne \emptyset$, and for any two positive real numbers $s, t$. 
Thanks to Lemma \ref{jacobian}, we know that the Jacobian ideal $J$ is contained in $I$. Compared with $J$, how big is this ideal $I$? Furthermore, is it possible to interpret $I$ geometrically? 
\end{rem}

Furthermore, thanks to Theorem \ref{multiplier}, Theorem \ref{subtau} gives a subadditivity formula for multiplier ideals on singular varieties. 
\begin{thm}\label{multsub}
Let $X$ be a normal $\Q$-Gorenstein variety over a field $K$ of characteristic zero. 
Let $\mathfrak{J}(X/K)$ be the Jacobian ideal sheaf of $X$ over $K$. 
Then, for any two nonzero ideal sheaves $\a, \b$ of $X$, and for any two positive real numbers $s, t$, 
$$\mathfrak{J}(X/K)\J(\a^t\b^s) \subseteq \J(\a^t)\J(\b^s). $$
\end{thm}

\begin{eg}
Consider the $\mathrm{A}_{2n}$-singularity $R=\C[X,Y,Z]/(XY-Z^{2n+1})$ with $n \ge 2$.
The Jacobian ideal $\mathfrak{J}(R/\C)$ is $(x,y,z^{2n})$. 
Let $\a=(x, y^{2n}, y^{2n-1}z, \dots, z^{2n})$. 
Then $\J(\a)=\a=(x, y^{2n}, y^{2n-1}z, \dots, z^{2n})$ and $\J(\a^{1/2})=(x,y^n, y^{n-1}z, \dots, z^n)$. 
Therefore, $\mathfrak{J}(R/\C)$ multiplies $\J(\a)$ into $\J(\a^{1/2})^2$ but $(x,y,z)$ does not because $xz \notin \J(\a^{1/2})^2$. 
Thus we cannot replace the Jacobian ideal $\mathfrak{J}(R/K)$ by its radical ideal $\sqrt{\mathfrak{J}(R/K)}$ in Theorem \ref{multsub}.
\end{eg}

\section{Summation formula on singular varieties}
When a ring is not necessarily regular, we prove a summation property of the ideal $\widetilde{\tau}((\a+\b)^t)$ with the aid of Lemma \ref{key lemma}.
\begin{thm}\label{tausum}
Let $R$ be an F-finite reduced ring of characteristic $p>0$ and $\a, \b \subseteq R$ be two ideals such that $\a \cap R^{\circ} \ne \emptyset$ and $\b \cap R^{\circ} \ne \emptyset$. 
Fix any real number $t>0$. Then
$$\widetilde{\tau}((\a+\b)^t) = \sum_{\lambda+\mu=t}\widetilde{\tau}(\a^{\lambda}\b^{\mu}).$$
\end{thm}

\begin{proof}
By definition, it immediately follows that $\widetilde{\tau}((\a+\b)^t) \supset \widetilde{\tau}(\a^{\lambda}\b^{\mu})$ for any $\lambda, \mu \ge 0$ with $\lambda+\mu=t$, because $\a^{\lceil \lambda p^e \rceil}\b^{\lceil \mu p^e \rceil} \subset (\a+\b)^{\lceil tp^e \rceil}$ for every $q=p^e$. 


Therefore, we will show $\widetilde{\tau}((\a+\b)^t) \subseteq \sum_{\lambda+\mu=t}\widetilde{\tau}(\a^{\lambda}\b^{\mu})$.  
We may assume without loss of generality that $R$ is a local ring. 
Fix an arbitrary element $c \in \widetilde{\tau}((\a+\b)^t)$. 
Let $x^{(e)}_1, \dots, x^{(e)}_{r_e}$ be a system of generators of $(\a+\b)^{\lceil tq \rceil}$ for each $q=p^e$. 
We may assume that $x^{(e)}_i$ belongs to an ideal $\a^{k^{(e)}_i}\b^{l^{(e)}_i}$ for some integers $k^{(e)}_i, l^{(e)}_i \ge 0$ such that $k^{(e)}_i+l^{(e)}_i=\lceil tp^e \rceil$ for all $e \ge 0$ and all $1 \le i \le r_e$. 
Let $d \in R^{\circ}$ be an element such that the localized ring $R_d$ is regular. 
Then, by Theorem \ref{exist} (1), some power $d^n$ of $d$ is an element of $\widetilde{\tau}(R)$. 
By Lemma \ref{key lemma}, there exist an integer $e'>0$ and $R$-linear maps $\phi^{(e)}_i \in \Hom_R(R^{1/p^e},R)$ for $0 \le e \le e'$ and $1 \le i \le r_e$ such that 

$$c=\sum_{e=0}^{e'}\sum_{i=1}^{r_e} \phi^{(e)}_i((d x^{(e)}_i)^{1/p^e}).$$
Denote $\lambda^{(e)}_i=k^{(e)}_i/p^e$ and $\mu^{(e)}_i=(l^{(e)}_i+tp^e-\lceil tp^e \rceil)/p^e$ for all $0 \le e \le e'$ and all $1 \le i \le r_e$. 
Since $\lceil \lambda^{(e)}_i p^e \rceil=k^{(e)}_i$, $\lceil \mu^{(e)}_i p^e \rceil=l^{(e)}_i$ and $\lambda^{(e)}_i+\mu^{(e)}_i=t$, by Lemma \ref{key lemma} again, 
$$\phi^{(e)}_i((d^nx^{(e)}_i)^{1/p^e}) \in \widetilde{\tau}(\a^{\lambda^{(e)}_i}\b^{\mu^{(e)}_i}) \subset \sum_{\lambda+\mu=t}\widetilde{\tau}(\a^{\lambda}\b^{\mu})$$
for all $0 \le e \le e'$ and all $1 \le i \le r_e$. 
Thus we have $c \in \sum_{\lambda+\mu=t}\widetilde{\tau}(\a^{\lambda}\b^{\mu})$. 
\end{proof}

By virtue of Theorem \ref{multiplier}, the above formula for the ideal $\widetilde{\tau}((\a+\b)^t)$ leads us to the following generalization of Musta\c t\v a's summation formula \cite{Mu3} for multiplier ideals. 

\begin{thm}\label{multsum}
Let $X$ be a normal $\Q$-Gorenstein variety over a field $K$ of characteristic zero and $\a, \b \subseteq \O_X$ be two nonzero ideal sheaves on $X$. Fix any real number $t > 0$. Then
$$\J(X, (\a+\b)^t)=\sum_{\lambda+\mu=t}\J(X,\a^{\lambda}\b^{\mu}).$$
In particular, by Theorem \ref{multsub}, 
$$\mathfrak{J}(X/K) \J((\a+\b)^t) \subseteq \sum_{\lambda+\mu=t}\J(\a^{\lambda}) \J(\b^{\mu}),$$
where $\mathfrak{J}(X/K)$ is the Jacobian ideal of $X$ over $K$. 
\end{thm}

\begin{proof}
First we will show that $\J(X, (\a+\b)^t) \supset \J(X,\a^{\lambda}\b^{\mu})$ for all real numbers $\lambda, \mu \ge 0$ with $\lambda+\mu=t$. 
Take a common log resolution $f:\widetilde{X} \to X$ of $\a, \b$ and $\a+\b$ so that $\a\O_{\widetilde{X}}=\O_{\widetilde{X}}(-F_a)$, $\b\O_{\widetilde{X}}=\O_{\widetilde{X}}(-F_b)$ and $(\a+\b)\O_{\widetilde{X}}=\O_{\widetilde{X}}(-F_c)$ are invertible, and write $F_a=\sum a_iE_i$, $F_b=\sum b_iE_i$ and $F_c=\sum c_iE_i$. 
Then $c_i=\min\{a_i, b_i\}$. 
Thus, for all $\lambda, \mu \ge 0$ with $\lambda+\mu=t$, one has 
$$\lceil K_{\widetilde{X}/X} - tF_c \rceil \ge \lceil K_{\widetilde{X}/X}-\lambda F_a - \mu F_b \rceil, $$
which implies the inclusion $\J(X, (\a+\b)^t) \supset \J(X,\a^{\lambda}\b^{\mu})$. 

Now we will prove the reverse inclusion 
$\J(X, (\a+\b)^t) \subseteq \sum_{\lambda+\mu=t}\J(X,\a^{\lambda}\b^{\mu})$. 
The question is local, so we may assume that $X=\Spec R$, where $R$ is a $\Q$-Gorenstein normal local ring of essentially of finite type over a field of characteristic zero. 
We denote by $(\widetilde{R}, \widetilde{\a}, \widetilde{\b})$ reduction to characteristic $p \gg 0$ of $(R, \a, \b)$ (see Setup \ref{reduction}). 
Then it follows from Theorem \ref{tausum} and \cite[Theorem 6.8, 6.9]{HY} that 
\begin{align*}
\J(\Spec \widetilde{R}, (\widetilde{\a}+\widetilde{\b})^t)=\tau(\widetilde{R}, (\widetilde{\a}+\widetilde{\b})^t) & = \sum_{\lambda+\mu=t}\tau(\widetilde{R},\widetilde{\a}^{\lambda}\widetilde{\b}^{\mu})\\
& \subseteq \sum_{\lambda+\mu=t}\J(\Spec \widetilde{R}, \widetilde{\a}^{\lambda}\widetilde{\b}^{\mu}).
\end{align*}
This means a similar inclusion holds in characteristic zero, that is, 
$$\J(\Spec R, (\a+\b)^t) \subseteq \sum_{\lambda+\mu=t}\J(\Spec R,\a^{\lambda}\b^{\mu}).$$ 
\end{proof}


\section{Asymptotic case and its application}
Ein, Lazarsfeld and Smith \cite{ELS} introduced the notion of asymptotic multiplier ideals, which is a variant of multiplier ideals defined for a graded family of ideals. 
In this paper, a graded family $\a_{\bullet}=\{ \a_m \}_{m \ge 1}$ of ideals on a Noetherian ring $R$ means a collection of ideals $\a_m \subseteq R$, satisfying $\a_1 \cap R^{\circ} \ne \emptyset$ and $\a_k \cdot \a_l \subseteq \a_{k+l}$ for all $k, l \ge 1$. 
Just for convenience, we decree that $\a_0=R$. 
A graded family of ideals on an algebraic variety is also defined similarly. 
One of the most important examples of graded families of ideals is a collection of symbolic powers $\a^{(\bullet)}=\{\a^{(m)}\}_{m \ge 1}$.

Let $\a_{\bullet}=\{\a_m\}$ be a graded family of ideals on a $\Q$-Gorenstein normal variety $X$ over a field of characteristic zero and let $t>0$ be a real number. 
Then the asymptotic multiplier ideal sheaf $\J(\a_\bullet^t)$ of $\a_\bullet$ with exponent $t$ is defined to be the unique maximal member among the family of ideals $\{\J(\a_m^{t/m})\}$ for $m \ge 1$ with respect to inclusion. 
The reader is referred to \cite{ELS} and \cite{La} for details (cf. Remark \ref{maximal}). 
Let $\b_{\bullet}$ be another graded family of ideals on $X$ and let $s>0$ be another real number. 
Then the ``mixed'' asymptotic multiplier ideal sheaf $\J(\a_\bullet^t \b_{\bullet}^s)$ is defined to be the unique maximal member among the family of ideals $\{\J(\a_m^{c/m}\b_m^{d/m})\}$ for $m \ge 1$. 

Recently Hara defined a variant of the ideal $\tau(\a^t)$ corresponding asymptotic multiplier ideals. In order to define this variant, Hara introduced a variant of $\a^t$-tight closure defined for a graded family of ideals. 

\begin{defn}[\textup{\cite[Definition 2.7]{Ha}}]
Let $\a_{\bullet}=\{a_m \}$ be a graded family of ideals on a Noetherian reduced ring $R$ of characteristic $p>0$ and let $t>0$ be a real number. 
Let $N \subseteq M$ be $R$-modules. 
The $\a_{\bullet}^t$-tight closure of $N$ in $M$, denoted by $N_M^{*\a_{\bullet}^t}$, is defined to be the submodule of $M$ consisting of all elements $z \in M$ for which there exists $c \in R^{\circ}$ such that 
$$c \a_{\lceil tq \rceil} z^q \subseteq N_M^{[q]}$$
for all large $q=p^e$. 
\end{defn} 

\begin{rem}[\textup{\cite[Observation 2.8]{Ha}}]\label{maximal}
Let $\a_{\bullet}=\{a_m\}$ be a graded family of ideals on Noetherian ring $R$ of characteristic $p>0$. Then given any real number $t>0$ one has the inclusion $\tau(\a_k^{t/k}) \subseteq \tau(\a_{kl}^{t/kl})$ (resp. $\widetilde{\tau}(\a_k^{t/k}) \subseteq \widetilde{\tau}(\a_{kl}^{t/kl})$) for all integers $k,l \ge 1$ (because $c\a_k^{\lceil tq/k \rceil} \subseteq \a_{kl}^{\lceil tq/kl \rceil}$ for some $c \in R^{\circ}$ and every $q=p^e$). 
This implies that the family of ideals $\{\tau(\a_m^{t/m})\}_{m \ge 1}$ (resp. $\{\widetilde{\tau}(\a_m^{t/m})\}_{m \ge 1}$) has a unique maximal element with respect to inclusion: the existence of at least one maximal member follows from the ascending chain condition on ideals (because $R$ is Noetherian). 
\end{rem}

\begin{propdef}[\textup{\cite[Proposition-Definition 2.9]{Ha}}]
Let $\a_{\bullet}=\{a_m \}$ be a graded family of ideals on a Noetherian excellent reduced ring $R$ of characteristic $p>0$ and let $t>0$ be a real number. 
Let $E =\bigoplus_{\m} E_R(R/\m)$ be the direct sum, taken over all maximal ideals $\m$ of $R$, of the injective hulls of the residue fields $R/\m$. 
We define 
$$\tau(\a_{\bullet}^t)=\tau(t \cdot \a_{\bullet})=\bigcap_{M \subset E} \Ann_R (0_M^{*\a_{\bullet}^t}),$$
where the intersection on the right hand side is taken over all finitely generated $R$-submodules $M$ of $E$.  Also we define 
$$\widetilde{\tau}(\a_{\bullet}^t)= \widetilde{\tau}(t \cdot \a_{\bullet})=\Ann_R (0_E^{*\a_{\bullet}^t}).$$
Then $\tau(\a_{\bullet}^t)$ $($resp. $\widetilde{\tau}(\a_{\bullet}^t)$ $)$ is equal to the unique maximal element among the set of ideals $\{\tau((\a_m)^{t/m}\}_{m \ge 1}$ $($resp. $\{\widetilde{\tau}((\a_m)^{t/m}\}_{m \ge 1}$ $)$ with respect to inclusion. 
\end{propdef}

Then we have a formula for ``asymptotic $\tau$'' similar to Theorem \ref{subtau}. 
\begin{prop}\label{symblem}
Let $R$ be an equidimensional reduced affine algebra over a perfect field $K$ of characteristic $p>0$ and let $J=\mathfrak{J}(R/K)$ be the Jacobian ideal of $R$ over $K$. 
Let $\a_{\bullet}=\{\a_m\}$ be a graded family of ideals on $R$.
Fix positive integers $k$ and $l$, and a real number $t>0$. 
Then 
$$J \widetilde{\tau}(\a_{\bullet}^{t(k+l)}) \subseteq  \widetilde{\tau}(\a_{\bullet}^{tk})\widetilde{\tau}(\a_{\bullet}^{tl}).$$
In particular, 
$$J^{l-1} \widetilde{\tau}(\a_{\bullet}^{tkl}) \subseteq \widetilde{\tau}(\a_{\bullet}^{tk})^l$$
for every $l>0$. 
\end{prop}
\begin{proof}
We can choose sufficiently large $m\gg 0$ so that $\widetilde{\tau}(\a_{\bullet}^{t(k+l)})=\widetilde{\tau}(\a_{m}^{t(k+l)/m})$. 
Then, by our subadditivity formula (Theorem \ref{subtau}), we have 
$$J \widetilde{\tau}(\a_{m}^{t(k+l)/m}) \subseteq \widetilde{\tau}(\a_{m}^{tk/m}) \widetilde{\tau}(\a_m^{tl/m}).$$
Thus $J \widetilde{\tau}(\a_{\bullet}^{t(k+l)}) \subseteq  \widetilde{\tau}(\a_{\bullet}^{tk})\widetilde{\tau}(\a_{\bullet}^{tl})$ as required. 
\end{proof}

Before stating an application of Proposition \ref{symblem}, we include the next lemma for further reference, whose proof is immediate from definition.   
\begin{lem}\label{symblem2}
Let $\a_{\bullet}=\{a_m \}$ be a graded family of ideals on a Noetherian reduced ring $R$ of characteristic $p>0$ and let $t>0$ be a real number. Then for all integers $k, l \ge 0$, 
$$\a_k \widetilde{\tau}(\a_{\bullet}^l) \subseteq \widetilde{\tau}(\a_{\bullet}^{k+l}).$$
\end{lem}
\begin{proof}
It is enough to show that $(0_E^{*\a_{\bullet}^l}:\a_k) \supseteq 0_E^{*\a_{\bullet}^{(k+l)}}$, but it is immediate because $\a_k^q \a_{lq} \subseteq \a_{(k+l)q}$ for every $q=p^e$. 
\end{proof}

As an application of our subadditivity formula, we obtain an answer to Hochster-Huneke's question \cite{HH5} concerning the growth of symbolic powers of ideals in singular affine algebras. 
\begin{thm}\label{symbthm}
Let $R$ be an equidimensional reduced affine algebra over a perfect field $K$ of characteristic $p>0$ and $J=\mathfrak{J}(R/K)$ be the Jacobian ideal of $R$ over $K$. 
Let $\a \subseteq R$ be any ideal such that $\a \cap R^{\circ} \ne \emptyset$.
Let $h$ be the largest analytic spread of $\a R_P$ as $P$ runs through the associated primes of $\a$. 
Then, for every integer $m \ge 0$ and every integer $n \ge 1$, 
$$\widetilde{\tau}(R)J^{n-1}\a^{(hn+mn)} \subseteq (\a^{(m+1)})^n. $$
In particular, one has
$$J^{n} \a^{(hn+mn)} \subseteq (\a^{(m+1))})^n$$
for all $n \ge 1$. 
\end{thm}

\begin{proof}
We consider the graded family $\a^{(\bullet)}=\{\a^{(m)}\}$ of symbolic powers of $\a$. 
By Proposition \ref{symblem} and Lemma \ref{symblem2}, 
$$\widetilde{\tau}(R)J^{n-1}\a^{(hn+mn)} \subseteq J^{n-1}\widetilde{\tau}((hn+mn) \cdot \a^{(\bullet)}) \subseteq \widetilde{\tau}((h+m) \cdot \a^{(\bullet)})^n.$$
Therefore, it suffices to show that $\widetilde{\tau}((h+m) \cdot \a^{(\bullet)}) \subseteq  \a^{(m+1)}$. 
Let $t$ be an indeterminate over $K$, let $L=K(t)$ and let $R_L=R \otimes_K L$. 
Note that the issues are unaffected by replacing $R$ by $R_L$ as in the proof of Lemma \ref{jacobian}. 
Thus we may assume without loss of generality that the residue field of each of the rings $R_P$ is infinite when $P$ is an associated prime of $\a$, and it follows that for each associated prime $P$ of $\a$, $\a R_P$ has a reduction ideal generated by at most $h$ elements. 
Then, by virtue of Proposition \ref{loc} and Theorem \ref{skoda}, 
\begin{align*}
\widetilde{\tau}((h+m) \cdot \a^{(\bullet)})R_P=\widetilde{\tau}((h+m) \cdot \a^{(\bullet)}R_P)&=\widetilde{\tau}((\a R_P)^{h+m})\\
&\subseteq \a^{m+1}R_P
\end{align*}
for every associated prime $P$ of $\a$, because after localization at $P$ the symbolic and ordinary powers of $\a$ are the same.  
Thus one has $\widetilde{\tau}((h+m) \cdot \a^{(\bullet)}) \subseteq \a^{(m+1)}$, as required. 

For the latter assertion, it is enough to show that $J$ is contained in $\widetilde{\tau}(R)$.  
Although it follows from essentially the same argument as the proof of \cite[(1.5.5)]{HH5} which states that $J \subseteq \tau(R)$, we give a slightly different proof here with the aid of Lemma \ref{key lemma}. 
Fix any element $c \in J$. 
By Theorem \ref{jacobianlem} and an argument similar to the proof of Lemma \ref{jacobian}, we may assume without loss of generality that there exists a regular subring $A$ of $R$ such that $R$ is module-finite and generically \'etale over $A$, $c R^{1/q} \subseteq A^{1/q}[R]$ and $A^{1/q}[R] \cong A^{1/q} \otimes_A R$ is flat over $R$ for every $q=p^e>0$. 
Moreover, we can assume that $R$ is local. 
Let $d \in R^{\circ}$ be an element of $A \cap \widetilde{\tau}(R)$. 
Then there exists a power $q$ of $p$ such that $R \xrightarrow{d^{1/q}} A^{1/q}[R]$ splits as an $R$-module homomorphism (because $A$ is regular: cf. \cite{HH2}), so we have an $R$-linear map $R^{1/q} \to R$ sending $d^{1/q}$ to $c$. 
By Lemma \ref{key lemma}, this implies that $c \in \widetilde{\tau}(R)$. 
\end{proof}

\begin{rem}
Hochster and Huneke illustrate in \cite[Example 3.8]{HH4} that the exponent $n-1$ used on the Jacobian ideal in Theorem \ref{symbthm} cannot be replaced by $n-2$ in general: the exponent $n-1$ in Theorem \ref{symbthm} is best possible. 
\end{rem}

Via reduction to characteristic $p \gg 0$, we obtain similar uniform bounds on the behavior of symbolic powers of ideals in affine algebras of characteristic zero. 
\begin{thm}
Let $R$ be a normal $\Q$-Gorenstein affine domain over a field $K$ of characteristic zero and $J=\mathfrak{J}(R/K)$ be the Jacobian ideal of $R$ over $K$. 
Let $\a \subseteq R$ be any nonzero ideal.
Let $h$ be the largest analytic spread of $\a R_P$ as $P$ runs through the associated primes of $\a$. 
Then, for every integer $m \ge 0$ and every integer $n \ge 1$, 
$$\J(R)J^{n-1}\a^{(hn+mn)} \subseteq (\a^{(m+1)})^n. $$
In particular, one has 
$$J^{n}\a^{(hn+mn)} \subseteq (\a^{(m+1)})^n$$
for all $n \ge 1$. 
\end{thm}

\vspace*{1em}
In the latter half of this section, we study a summation property of asymptotic multiplier ideals (resp. ``asymptotic $\tau$'') analogous to Theorem \ref{multsum} (resp. Theorem \ref{tausum}).

Given two graded families $\a_{\bullet}=\{a_m\}$ and $\b_{\bullet}=\{b_m\}$, the graded family $\a_{\bullet}+\b_{\bullet}=\{(\a_{\bullet}+\b_{\bullet})_m\}$ is given by
$$(\a_{\bullet}+\b_{\bullet})_m=\sum_{k+l=m}\a_k\b_l.$$

\begin{prop}\label{asymtausum}
Let $\a_{\bullet}=\{a_m \}$ and $\b_{\bullet}=\{b_m \}$ be graded families of ideals on an F-finite reduced ring $R$ of characteristic $p>0$. For each real number $t>0$, one has 
$$\widetilde{\tau}((\a_{\bullet}+\b_{\bullet})^t)=\sum_{\lambda+\mu=t}\widetilde{\tau}(\a_{\bullet}^\lambda \b_{\bullet}^\mu).$$
\end{prop}
\begin{proof}
First we will show that $\widetilde{\tau}((\a_{\bullet}+\b_{\bullet})^t) \supset \widetilde{\tau}(\a_{\bullet}^\lambda \b_{\bullet}^\mu)$ for any $\lambda, \mu \ge 0$ with $\lambda+\mu=t$. 
Take sufficiently  large $m \gg 0$ so that 
$\widetilde{\tau}(\a_{\bullet}^\lambda \b_{\bullet}^\mu)=\widetilde{\tau}(\a_m^{\lambda/m}\b_m^{\mu/m})$. 
Then, by Theorem \ref{tausum}, 
$$\widetilde{\tau}(\a_m^{\lambda/m}\b_m^{\mu/m}) \subset \widetilde{\tau}((\a_m+\b_m)^{t/m}) \subset \widetilde{\tau}((\sum_{k+l=m}\a_k \b_l)^{t/m}).$$
Thus we have $\widetilde{\tau}(\a_{\bullet}^\lambda \b_{\bullet}^\mu) \subset \widetilde{\tau}((\a_{\bullet}+\b_{\bullet})^t)$. 

Conversely, we will prove that $\widetilde{\tau}((\a_{\bullet}+\b_{\bullet})^t) \subseteq \sum_{\lambda+\mu=t}\widetilde{\tau}(\a_{\bullet}^\lambda \b_{\bullet}^\mu)$. 
If $\m$ is sufficiently large and divisible, then by Theorem \ref{tausum} again, 
$$\widetilde{\tau}((\a_{\bullet}+\b_{\bullet})^t)=\widetilde{\tau}((\sum_{k+l=m}\a_k\b_l)^{t/m})=\sum_{t_0+ \cdots t_m=t/m}\widetilde{\tau}(\prod_{i=0}^m (\a_i\b_{m-i})^{t_i}).$$
Since $(\a_i)^{m!/i} \subseteq \a_{m!}$ for all $1 \le i \le m$ and $(\b_{m-i})^{m!/(m-i)} \subseteq \b_{m!}$ for all $0 \le i \le m-1$, 
$$\widetilde{\tau}(\prod_{i=0}^m (\a_i\b_{m-i})^{t_i}) \subseteq \widetilde{\tau}(\prod_{i=0}^m \a_{m!}^{\frac{it_i}{m!}}\b_{m!}^{\frac{(m-i)t_i}{m!}})=\widetilde{\tau}(\a_{m!}^{\sum_{i=0}^m \frac{it_i}{m!}} \b_{m!}^{\sum_{i=0}^m \frac{(m-i)t_i}{m!}}).$$ 
Since $\sum_{i=0}^m \frac{it_i}{m!}+\sum_{i=0}^m \frac{(m-i)t_i}{m!}=\frac{t}{m!}$ for all $t_0, \dots, t_m \ge 0$ with $\sum_{i=0}^m t_i=t$, one has 
$$\widetilde{\tau}((\a_{\bullet}+\b_{\bullet})^t) \subseteq \sum_{\lambda+\mu=t}\widetilde{\tau}(\a_{m!}^{\lambda /m!} \b_{m!}^{\mu /m!}) \subseteq \sum_{\lambda+\mu=t}\widetilde{\tau}(\a_{\bullet}^{\lambda} \b_{\bullet}^{\mu}).$$
\end{proof}

We can easily derive from Theorem \ref{multsum} a generalization of Musta\c t\v a's summation formula for asymptotic multiplier ideals, while the original proof of his formula is a little complicated. 
\begin{prop}
Let $\a_{\bullet}=\{a_m \}$ and $\b_{\bullet}=\{b_m \}$ be graded families of ideals on a $\Q$-Gorenstein normal variety $X$ over a field of characteristic zero. For each real number $t>0$, one has 
$$\J((\a_{\bullet}+\b_{\bullet})^t)=\sum_{\lambda+\mu=t}\J({\a_{\bullet}}^\lambda{\b_{\bullet}}^\mu).$$
In particular, 
$$\mathfrak{J}(X/K) \J((\a_{\bullet}+\b_{\bullet})^t) \subseteq \sum_{\lambda+\mu=t}\J({\a_{\bullet}}^\lambda)\J({\b_{\bullet}}^\mu),$$
where $\mathfrak{J}(X/K)$ is the Jacobian ideal sheaf of $X$ over $K$. 
\end{prop}

\begin{proof}
The proof is essentially the same as that of Proposition \ref{asymtausum}. 
\end{proof}

\begin{rem}
Statements similar to the results in Section 2 and 3 also hold for ``asymptotic $\tau$," 
and these give alternative proofs of the results in Section 4.   
\end{rem}

\begin{acknowledgement}
The author is indebted to Melvin Hochster for helpful advice and Robert Lazarsfeld for warm encouragement on this work. 
He is also grateful to Ken-ichi Yoshida for useful comments and valuable conversations. 
A part of this work was done at the University of Michigan where the author stayed in the fall of the year 2004. 
The author would like to express his deep gratitude to the University of Michigan for the hospitality and support. 

\end{acknowledgement}

\end{document}